**Michiel Hazewinkel**   1   CWI
Direct line: +31-20-5924204   POBox 94079
Secretary: +31-20-5924233   1090GB Amsterdam
Fax: +31-20-5924166
E-mail: mich@cwi.nl   original version: 6 November, 1994
revised version: 24 oktober 2004

# Hopf algebras: their status and pervasiveness
## (as of Oct. 2004)

by

*Michiel Hazewinkel*
*CWI*
*POBox 94079*
*1090GB Amsterdam*
*The Netherlands*

In late October 2004 I did a search on Hopf algebras in the database ZMATH to find out where they occur and what use is made of them (I.e. what applications there are). The result is a little astonishing, as can be seen from what follows.

First of all, here are the numbers.

Overall: 10589 hits. If the search is limited to the period up to and including 1994, the score is 3828. I.e. about two thirds of Hopf algebra related papers appeared in the last 10 years. Similar results appear when the preprint archive 'arXive' is consulted.

These numbers are not very precise. There is in fact quite a good deal more: it is perfectly possible for an (abstract of an) article to be about Hopf algebras or quantum groups without ever mentioning any of those phrases exactly. For instance there are 18 papers about quantum polynomials that are not among the 10589 mentioned above and 742 papers on quantum spaces that are not among these.

Below are the number of hits for each of the main subfields of mathematics (where the very general subfields 00: general; 01: History and biography; 97: Education have been left out.

The search was done on "('Hopf algebras' or 'bialgebras' or 'quantum groups' or 'quantum algebras') and 'main classification number'". 'Quantum groups' are a special kind of Hopf algebras; 'bialgebras' are slightly more general (and many bialgebras are automatically Hopf algebras), 'quantum algebras' are in principle more general than bialgebras and Hopf algebras; however, including or excluding them from the search makes a difference of only about 4%.
    The phrases 'enveloping algebras' and 'group algebras' (group rings) were not included in the search. Both are special kinds of Hopf algebras. But the special extra structure that makes them Hopf algebras is ignored in the majority of papers dealing with these objects. And, if not, the phrase 'Hopf algebra(s)' is likely to turn up, in the basic index, or the classification 16W30 will appear.

| | |
|---|---:|
| 03 Mathematical logic and foundations | 23 |
| 05 Combinatorics | 282 |
| 06 Order, lattices, ordered algebraic structure | 38 |
| 08 General algebraic systems | 18 |
| 11 Number theory | 122 |
| 12 Field theory and polynomials | 47 |
| 13 Commutative rings and algebras | 123 |
| 14 Algebraic geometry | 553 |
| 15 Linear and multilinear algebra. Matrix theory | 126 |
| 16 Associative rings and algebras | 3699 |
| 17 Nonassociative rings and algebras | 5199 |

| | | |
|---|---|---:|
| 18 | Category theory; abstract homological algebra | 589 |
| 19 | K-theory | 103 |
| 20 | Groups | 710 |
| 22 | Topological groups, Lie groups | 426 |
| 26 | Real functions | 18 |
| 28 | Measure and integration | 6 |
| 30 | Functions of a complex variable | 12 |
| 31 | Potential theory | 0 |
| 32 | Several complex variables and analytic spaces | 89 |
| 33 | Special fucntions | 389 |
| 34 | Ordinary differential equations | 25 |
| 35 | Partial differential equations | 187 |
| 37 | Dynamical systems and ergodic theory | 338 |
| 39 | Difference and functional equations | 110 |
| 40 | Sequences, series, summability | 1 |
| 41 | Approximations and expansions | 7 |
| 42 | Fourier analysis | 12 |
| 43 | Abstract harmonic analysis | 73 |
| 44 | Integral transforms, operational calculus | 5 |
| 45 | Integral equations | 1 |
| 46 | Functional analysis | 989 |
| 47 | Operator theory | 122 |
| 49 | Calculus of variations and optimal control; optimization | 2 |
| 51 | Geometry | 11 |
| 52 | Convexity | 20 |
| 53 | Differential geometry | 427 |
| 54 | General topology | 13 |
| 55 | Algebraic topology | 411 |
| 57 | Manifolds and cell complexes | 792 |
| 58 | Global analysis, analysis on manifolds | 598 |
| 60 | Probability theory and stochastic processes | 104 |
| 62 | Statistics | 3 |
| 65 | Numerical analysis | 14 |
| 68 | Computer science | 35 |
| 70 | Mechanics of particles and systems | 40 |
| 74 | Mechanics of deformable solids | 4 |
| 76 | Fluid mechanics | 21 |
| 78 | Optics, electromagnetic theory | 8 |
| 80 | Classical theromodynamics, heat transfer | 3 |
| 81 | Quantum theory | 5463 |
| 82 | Statistical mechanics, structure of matter | 833 |
| 83 | Relativity and gravitational theory | 194 |
| 85 | Astronomy and astrophysics | 2 |
| 86 | Geophysics | 0 |
| 90 | Operations research, mathematical programming | 1 |
| 91 | Game theory, economics, social and behavioural sciences | 3 |
| 92 | Biology and other natural sciences | 9 |
| 93 | Systems theory; control | 8 |
| 94 | Information and communication, circuits | 9 |

Here are some words on what kind of Hopf algebra applications the numbers above allude to.

    * The exact phrase "Hopf algebras" occurs twice in the MSCS: 16W30 and 57T05. Basically a Hopf agebra is an algebraic structure; more precisely it is an algebra with additional structure. This gives the classification 16W30 (= coalgebras, bialgebras, Hopf algebras and modules on which they act). The reason that there is a classification 57T05: Hopf algebras, is

that historically the first examples of Hopf algebras came in the form of the cohomology or homology of suitable spaces and manifolds.

* 03. Applications of Hopf algebras etc. to logic include applications of operads in categorical logic and use of the fact that categories of representations of Hopf algebras provide models for linear logic (noncommutative).

* 05 stands for the subfield ''Combinatorics'. A good many combinatorial identities arise by looking at representations of Hopf algebras in two ways (often universal enveloping algebras, a special kind of Hopf algebras), whence lots of applications of Hopf algebras in subsubfield 05A19 (= Combinatorial identities).

Many combinatorial objects naturally form the basis of a Hopf algebra; for instance planar binary trees and chord diagrams. Whence e.g. the role of Hopf algebras in 05C (= Graph theory).

The symmetric functions and several generalizations are best looked at as forming Hopf algebras. Hence the role of Hopf algebras in subfield 05E (Algebraic combinatorics).

There are also relations between Baxter algebras and Hopf algebras (notably free Baxter algebras).

* 06. As in the case of field 05 there are natural Hopf algebra structures with as basis various classes of posets. Other applications include incidence Hopf algebras and Galois type connections. Finally quantum logic relates to quantum algebras (just as there are algebras attached to classical logics).

* 08. Just as there are 3- and n-groups (sets with an operation that associates to each triple of elements a new element etc.) and n-Lie algebras there are corresponding quantum and Hopf structures. As far as I know this is just generalization for generalization's sake. More important are the links of bialgebras with duality contexts.

* 11, 12. There are quite a few applications of Hopf algebras etc. here.

First there the applications of formal groups (Lubin-Tate formal groups a.o.) to class field theory and reciprocity maps.

Second there are the applications of these formal groups and others (such as the fake monster formal group) to congruences of the coeffients of modular and automorphic forms and other aspects of automorphic forms such as the monstrous moonshine conjectures.

Third there is 'Hopf-Galois theory'. Classical Galois theory associates a group to certain field extensions and then intermediate fields correspond to subgroups. This does not work for inseparable extensions but at least in the case of height one purely inseparable extensions Lie algebras work. Hopf algebras generalize both groups and Lie algebras and so it is natural to attempt to study things in this contexts. This works and goes far beyond the combination of the two classical cases. An important aspect is that that there is usually more than one Hopf algebra that can do the job.

Finally there is a remarkable connection of Hopf algebras (free coalgebras really) to recursive sequences. See also 93.

* 13. Many important Hopf algebras are commutative. So one would expect a subfield 13??? (13 = Commutative rings and algebras), for commutative algebras with extra structure and a place for Hopf algebras there. But there isn't.

On the other hand there is an entry for the Witt vectors. The big Witt vectors are represented by *Symm*, the Hopf algebra of symmetric functions, arguably the most beautiful and rich object in present day mathematics. Immediately related are lambda and beta rings; for instance because the universal lambda ring on one generator is again *Symm*.

There is also work on an algebraic variant of Hodge decomposition in the context of cohomology Hopf algebras of algebras.

Still other aspects are Hopf-Galois extensions of commutative rings, Hopf algebras as free objects in categories of Baxter algebras and relations with differential and difference algebras.

* 14. Subsubfield 14L05 is 'Formal groups, p-divisible groups', a subsubfield of 14 (= Algebraic geometry). Formal groups are a special kind of Hopf algebras which in turn have applications in number theory (11M, 11R, 11S), field theory (12F). Another kind of Hopf

algebras is formed by (the coordinate rings of) algebraic groups (14L, 20G, 20H).

Quantum cohomolgy and Gromov Witten invariants constitute another application area here.

* 15. A major theme here is a quantum version of invariant theory (which also turns up in areas 14, 20, 22). Besides that there are quantum versions of Clifford algebras and papers on quantum random matrices.

* 16. Structure theory, properties, recognition theorems, reconstruction theorems (`a la Tannaka-Krein) etc. of Hopf algebras and bialgebras belong here (16W30). In addition to that (the coordinate rings of) quantum groups are a highy structured and relatively well understood class of noncommutative algebras. As such they serve as a valuable class of examples for testing a number of conjectures in noncommutative ring theory.

* 17B stands for 'Lie algebras', and 17B10 for the algebraic approach to representations of Lie algebras. A representation of a Lie algebras is the same thing as a representation of its universal enveloping algebra which is a special kind of Hopf algebra.

17B37 is ''Quantum groups, quantized enveloping algebras and related deformations'. This kind of defomations is also a special kind of Hopf algebra. Naturally there has been a lot of research trying to extend to quantum deformed enveloping algebras the kind of results known for classical enveloping algebras (resulting e.g. in a quantum PBW theorem). But this is not all by a long shot. As happens regularly when studying deformations of structures new insights appear. In this case crystal graphs and crystal bases. The quantum version of a Lie enveloping algebra depends on a parameter. When that parameter is 1 the classical situation is recovered. In a suitable formulation one can also set the parameter to zero. This leads to new insights (Kashiwara) for the undeformed Lie algebras (crystal bases). The zero and 'crystal' vaguely refer to physics at absolute zero.

* 18. The relevant words here seem to be multiple categories, braided categories, fusion categories, ribbon categories, etc. Also relations to noncommutative geometry. The papers recorded in this field seem (to me) to be preparations for possible future applications rather then real applications.

* 19. There is quantum K-theory. apart from that most papers listed here have to do with noncommutative geometry and cyclic and quantum cohomology.

* 20C is 'Representations of groups'. A representation of a group is the same thing as a representation of its group ring which is another special kind of Hopf algebra.

* 20G is the subfield of 'linear algebraic groups' and 20G42 is the subfield of 'quantum groups' which are special kinds of Hopf algebras (deformations of the coordinate rings of the classical groups).

* 22. This overlaps with 20 and 17. Apart from that there is the study of amenability of quantum groups.

* 26,28. Not much. Relations with p-adic analysis and nothing else.

* 30. q-analytic functions and such thing as q-versions of the Ward identities on Riemann surfaces.

* 32. The KZ equation and its quantum version as well as aspects of Mirror symmetry. For instance the matter of the monodromy of the KZ equation being the representation coresponding to the R-matrix of the relevant bracket algebra.

* 33. Orthogonal polynomials and other special functons have been around for something like 400 years. Already early in the early nineteenth century parameter dependent version of these were discovered called basic special functions or q-special functions, the parameter 'q' as was the usual one. This was three quarters of a century before 'quantum theory'. In the 1960's the seminal discovery was made by Vilenkin and W. Miller Jr (independently and practically simultaneously) that the classical special functions are essentially the coefficients of matrix representations of certain Lie groups. And this accounts for many of the remarkable properties of orthogonal poynomials (including that very orthogonality). In 1998-1989 it was discovered by three groups, independently and more or less at the same time, that the q-special fucntions relate to quantum groups in the same way as special functions relate to the classical Lie groups

(Vaksman, Soibelman, Ukraina; Masuda, Mimachi, Nakagami, Noumi, Ueno, Japan: Koornwinder, Netherlands). It is a nice coincidence that the 'q' used traditionally fits with the 'q' of quantum.

* 34, 35. Papers found here seem to be all about quantum completely integrable systems, the quantum inverse scattering method and related matters such as the q-version of the Riemann-Hilbert problem. In addition there are studies on quantum manifolds and qPDE's on them. For instance the quantum version of the KZ equation (Knizhnik-Zamolodchikov equation).

*39. q-versions of discrete versions of the things (such as completely integrable ssystems) that occur in 34, 35. Also quantum versions of the things that turn up when examining difference operators such as Hall-Littlewood functions.

* 40, 41. Nothing that is not recorded more appropriately elsewhere.

* 42. There is a way of looking at wavelets via a q-deformed algebra structure. Also of course there is a q-version of the Fourier transform. Also there is an intriguing way to use the quantum cohomology of Grassmannians to deal with interppolation problems for rational functions.

* 43. Here the key phrases appear to be: quantum Fourier transform, q-Lobachevskij space, quantum plane, q-analogs of the Riemann zeta function, stochastic processes on quantum groups, quantum homogeneous spcaces and symmetric domains (which really should be in 20, 22, and 32), Plancherel identity for Hopf algebras.

* 44. Just q-analogs of integral transforms

* 45. Just one case of a q-deformation of an integral equation.

* 46 is the specialism 'Functional analysis' with subfield 46L: 'Selfadjoint operator algebras'. There are well over 300 published papers dealing with Hopf algebras of this kind. Mostly in subsubfields 47L85 (Noncommutative topology), 46L87 (Noncommutative geometry). There are also some in 46M (Methods of category theory in functional analysis).

* 47. Nothing much. Quantum versions of some important things like exponentials. Also q-deformatons of operators and algebras of opeators and relations with Wick ordering and Wick algebras.

* 51. Just as groups relate to geometries (Erlangen program), quantum groups should relate to quantum geometries. There is a very little bit in this section about this. Basically this a large area that so far has not been explored.

* 52. There are Hopf algebras attached to polytopes.

* 53. This is a large area of applicatiions of quantum groups and Hopf algebras. It takes in quantum cohomology, Gromov-Witten invariants, noncommutative geometry, and differential calculi (different from the usual one) on manifolds.

There are also quantum manifolds which are modelled on quantum affine spaces in the same way that classical manifolds are modelled on the algebras of function on the standard vector spaces over the real or complex numbers and supermanifolds are modelled on algebras with both commuting and anticommuting variables.

* 54. Nothing here that is not (better) listed under 55.

* 55N is the subfield 'Homology and cohomology theories' of 55 (= Algebraic topology) and 55N22 is'Bordism, cobordism, formal group laws'. As already remarked 'a formal group is a special kind of Hopf algebra. For that matter the applications of formal group laws in (algebraic) topology are not limited to (co)bordism theories but also take in elliptic cohomology and (extraordinary) K-theories.

An important very sophisticated structure in algebraic topology ( a field that abounds in complicated tools and structures) is that of an operad. There are strong links between operads and Hopf algebras (bialgebras as free algebras over an operad).

Of course 55 is the field where Hopf algebras originally came from: as the (co)homology of suitable spaces (though the classification for that is now 57T05).

* 57M is 'Low dimensional topology'. There are numerous applications of Hopf algebras

(quantum groups) to 57M25 (= Knots and links) and also to other parts of 57M (such as invariants of three manifolds). More precisely link invariants are derived from R-matrices, that is solutions of the classical Yang-Baxter equation and these also define quantum groups. Conjecturally all link invariants arise (in one form or another) this way.

* 58 is the field ''Global analysis, analysis on manifolds', which has as subsubspecialisms 58M32 'Geometry of quantum groups' and 58M34 ''Noncommutative geometry (à la Connes)'. Here one also finds results on quantum homogeneous sapces and differential calculi associated to quantum groups.

* 60. Just as there are stochastic processes on (Lie) groups, there are qunatum stochastic processes on quantum groups.

I also believe that free probability ( D Voiculescu) has a Hopf algebra component. In any case since the work of Scheunert it has become clear that there is a strong link between Hopf algebas and coalgebras and stochastic processes.

* 68. Here Hopf algebras and bialgebras serve as carriers for considerations on semantics, (random) accessible data types, concurrent and distributed computing, and algebraic development techniques. Coalgebas are especially important for semantics.

* 70, 74, 76. Practically all that is listed under these classifications could be beter listed under 81 and 82. An exception might be the use of q-deformed Poisson structures in classical mechanics.

* 78. Nothing much except possibly the study of electromagnetism on q-spheres and spaces.

* 81. Quantum groups, quantum symmetry, and representations of quantum groups are a major topic here. The representation theory of the quantum deformations of the classical groups is very similar to that of the undeformed groups which goes some way towards explaining the large role of nonquantum Lie groups in quantum theory.

There is also the quantum Yang-Baxter equation and the quantum inverse scattering method (Bethe Ansatz).

Further, noncommutative geometry (again) and such things as Dirac operators on quantum manifolds.

Finally renormalization is controlled by the Connes-Kreimer Hopf algebra.

* 82. Here the relation with Hopf algebras lies again in R-matrices: exactly solvable lattice models. The Yang-Baxter equation derives part of its name from the work of R Baxter on exactly solvable models such as the eight vertex model.

* 85, 85. Besides matters that also occur under 81 and 82 there are here relations with quantum gravity. There is also a quantum Minkowsky space-time.

* 90. No applications or other relations except a reverse application of network flows to representation theory.

* 92. Nothing except use of crystal bases for the genetic code.

* 93. Some applications of Hopf algebra to realization theory. the essence of this lies probably in the coalgebra aspect in that free coalgebras exactly capture the recursiveness needed for finite dimensional realizations and recognizable power series (in several noncommuting variables) in the sense of Schützenberger.

* 94. Here some relations with Hopf algebras find expression through quantum information theory and quantum coding theory.

To conclude let me remark that the (more or less recently established) electronic preprint collection arXive has some 4230 preprints stored on Hopf algebras and quantum groups.